\magnification=\magstep1
\input amstex
\documentstyle{amsppt}

\define\defeq{\overset{\text{def}}\to=}
\define\ab{\operatorname{ab}}

\define\Gal{\operatorname{Gal}}

\def \isom {\overset \sim \to \rightarrow}
\define\Sect{\operatorname{Sect}}
\define\Spec{\operatorname{Spec}}

\def\cl{\operatorname{cl}}
\def\Card{\operatorname{Card}}
\def\Spf{\operatorname{Spf}}

\def\Fr{\operatorname{Fr}}
\def\sep{\operatorname{sep}}

\def\sing{\operatorname{sing}}

\NoRunningHeads
\NoBlackBoxes
\topmatter

\title
On sections of arithmetic fundamental groups of open $p$-adic annuli.
\endtitle

\author
Mohamed Sa\"\i di
\endauthor

\abstract
We show the {\it non-existence} of sections of arithmetic fundamental groups of {\it open $p$-adic annuli of small radii}, this implies the {\it non-existence}
of sections of arithmetic fundamental groups of {\it formal boundaries of formal germs of  $p$-adic curves}.
\endabstract

\endtopmatter

\document
\subhead
\S 0. Introduction/Statement of the Main Result
\endsubhead
Let $p\ge 2$ be a prime integer, $k$ a {\bf $p$-adic local field} (i.e., $k/\Bbb Q_p$ is a finite extension), with ring of integers $\Cal O_k$, uniformiser $\pi$, and residue field $F$.
Thus, $F$ is a finite field of characteristic $p$. 

Let $X\to \Spec \Cal O_k$ be a flat, proper, relative $\Cal O_k$-curve, with $X$ normal, 
and $X_k\defeq X\times _{\Spec \Cal O_k}\Spec k$ geometrically connected.
Assume $X(F)\neq \emptyset$. 
Let $x\in X^{\cl}(F)$ be a closed point, $\Cal O_{X,x}$ the local ring at $x$, $\widehat {\Cal O}_{X,x}$ its completion, and $E\defeq \widehat {\Cal O}_{X,x}
\otimes_{\Cal O_k}k= \widehat {\Cal O}_{X,x}[\frac{1}{\pi}]$. 
Write $\Cal X\defeq \Spec E$, which we assume to be geometrically connected. We shall refer to $\Cal X$ as the {\bf formal germ of $X$ at $x$}.

Let $\eta$ be a geometric point of $\Cal X$ with values in its generic point. Thus, $\eta$ determines an algebraic closure $\overline k$ of $k$, and a geometric point 
$\bar \eta$ of $\Cal X_{\overline k}\defeq \Cal X\times_{\Spec k}\Spec \overline k$. There exists a canonical exact sequence of profinite groups (cf. [Grothendieck], Expos\'e IX, Th\'eor\`eme 6.1)
$$1\to \pi_1(\Cal X_{\overline k},\bar \eta)\to \pi_1(\Cal X, \eta) \to G_k\to 1.\tag 1$$
Here, $\pi_1(\Cal X, \eta)$ denotes the arithmetic \'etale fundamental group of $\Cal X$ with base
point $\eta$, $\pi_1(\Cal X_{\overline k},\bar \eta)$ the \'etale fundamental group of $\Cal {X}_{\overline k}$ with base
point $\bar \eta$, and $G_k\defeq \Gal (\overline k/k)$ the absolute Galois group of $k$. 

The sequence (1) splits if $\Cal X(k)\neq \emptyset$. This is for example the case if the morphism $X\to \Spec \Cal O_k$ is smooth at $x$.
If $\Cal X(k)=\emptyset$; for instance if $X$ is stable and regular, and $x$ is an ordinary $F$-rational double point of $X_F\defeq X\times _{\Spec \Cal O_k}F$, the existence of sections $s:G_k\to 
 \pi_1(\Cal X, \eta)$ of the projection $\pi_1(\Cal X,\eta)\twoheadrightarrow G_k$ would provide examples of sections of the projection $\pi_1(X_k,\eta)\twoheadrightarrow G_k$
 which are {\it non-geometric} ($\eta$ induces a geometric point of $X_k$, denoted also $\eta$, via the morphism $\Cal X_k\to X_k$), 
 i.e., which do not arise from rational points. These in turn will provide {\it counter-examples to the $p$-adic version of the Grothendieck anabelian section conjecture}.
This prompt the following question.

\proclaim{Question A}
With the above notations, assume $\Cal X(k)=\emptyset$. Does the exact sequence (1) split?
\endproclaim

In this note we investigate the case where $\Cal X$ is a {\it $p$-adic open annulus}. Let 
$$A\defeq\Cal O_k[[S]],\ \ \ B\defeq A\otimes _{\Cal O_k}k=A \left [\frac {1}{\pi} \right],$$  
$D\defeq \Spf A$ is the formal standard open disc, and $\Cal D\defeq D_k=\Spec B$ its "generic fibre" which is the standard open disc centred at the point $"S=0"$. 
Let $\Bbb A^1_{k}=\Spec k[S]$, 
$Y_k=\Bbb P^1_k$ its smooth compactification with function field $k(S)$, and $Y=\Bbb P^1_{\Cal O_k}$ the smooth compactification of $\Bbb A^1_{\Cal O_k}=\Spec \Cal O_{k}[S]$. 
We shall identify $A$ with the completion of the local ring of $Y$ at the closed point $"S=0"$.
We have a natural morphism $\Cal D\to \Bbb P^1_k$, which induces an identification between the set of closed points of $\Cal D$ and the set
$$\{x\in \Bbb P^1_k : \vert S(x) \vert <1\}.$$

For an integer $n\ge 1$, let 
$$A_n\defeq \frac {\Cal O_k[[S,T]]}{(S^nT-\pi)},\ \ \  B_n\defeq A_n\otimes _{\Cal O_k}k,\ \  \text{and} \ \  \Cal C_n\defeq \Spec B_n.$$
The natural embedding $\Cal C_n\hookrightarrow \Cal D$ induces an identification between the set of closed points of $\Cal C_n$ and the open annulus
$$\{x\in \Cal D : \vert \pi \vert ^{\frac{1}{n}}<\vert S(x) \vert <1\}.$$ 

Further, let $P\defeq A_{(\pi)}$ be the localisation of $A$ at the ideal $(\pi)$, and $\widehat P$ the completion of $P$, which is a complete discrete valuation ring 
isomorphic to 
$$\Cal O_k[[S]]\{S^{-1}\} \defeq \Big\{\sum _{i=-\infty}^{\infty}a_iS^i : a_i\in \Cal O_k, \ \underset 
{i\to -\infty} \to  \lim \vert a_i \vert=0 \Big\},$$ 
where $\vert \ \vert$ is a normalised absolute value of $\Cal O_k$ (cf. [Bourbaki], $\S2$, 5).
Let $L\defeq \Fr (\widehat P)$ be the fraction field of $\widehat P$, and $\Cal C_{\infty}\defeq \Spec L$. We shall refer to $\Cal C_{\infty}$ as a {\bf formal boundary} of the formal germs $\Cal D$, and 
$\Cal C_i$ for $i\ge 1$.
We have natural scheme morphisms 
$$\Cal C_{\infty}\to \cdots \to \Cal C_{n+1}\to \Cal C_n\to \cdots \to \Cal C_1\to \Cal D\to \Bbb P^1_k.$$

Let $\eta$ be a geometric point of $\Cal C_{\infty}$, which induces a geometric point (denoted also $\eta$) of $\Cal C_n$ for $n\ge 1$. 
For $i\in \Bbb N\cup \{\infty\}$, we have an exact sequence of arithmetic fundamental groups
$$1\to \pi_1(\Cal C_{i,\overline k},\bar \eta)\to \pi_1(\Cal C_i, \eta) \to G_k\to 1,\tag 2$$
where $\pi_1(\Cal C_i, \eta)$ denotes the arithmetic \'etale fundamental group of $\Cal C_i$ with base
point $\eta$, $\pi_1(\Cal C_{i,\overline k},\bar \eta)$ the \'etale fundamental group of $\Cal {C}_{i,\overline k}\defeq C_i\times _{\Spec k}\Spec {\overline k}$ with base
point $\bar \eta$; which is induced by $\eta$, and $G_k\defeq \Gal (\overline k/k)$ the absolute Galois group of $k$. Here $\overline k$ is the algebraic closure of $k$ determined by $\eta$.

Our main result in this note is the following. 
\proclaim {Theorem A} We use notations as above. There exists an integer $N\ge 1$, such that for every integer $n\ge N$,
the projection $\pi_1(\Cal C_n,\eta)\twoheadrightarrow G_k$ {\bf doesn't split}. 
\endproclaim

The author ignores, for the time being, if the projection $\pi_1(\Cal C_1, \eta) \twoheadrightarrow G_k$ splits or not. 

As a corollary of Theorem A, we obtain the following.

\proclaim {Theorem B} The projection $\pi_1(\Cal C_{\infty},\eta)\twoheadrightarrow G_k$ {\bf doesn't split}. 
\endproclaim

One of the consequences of Theorems A, and  B, is that one can not produce examples of sections of hyperbolic curves over $p$-adic local fields, which arise from sections of arithmetic fundamental groups of boundaries of formal fibres, or open annuli with small radii. Those sections would be non-geometric, hence would provide counter-examples to the $p$-adic version of the Grothendieck anabelian section conjecture.

Finally we observe the following. For $i\in \Bbb N\cup \{\infty\}$, let 
$\pi_1(\Cal C_{i,\overline k},\bar \eta)^{\ab}$ be the maximal abelian quotient of $\pi_1(\Cal C_{i,\overline k},\bar \eta)$, and consider the push-out diagram
$$
\CD
1@>>> \pi_1(\Cal C_{i,\overline k},\bar \eta)@>>> \pi_1(\Cal C_i, \eta) @>>> G_k @>>> 1 \\
@. @VVV @VVV @| \\
1@>>> \pi_1(\Cal C_{i,\overline k},\bar \eta)^{\ab}@>>> \pi_1(\Cal C_i, \eta)^{(\ab)} @>>> G_k @>>> 1 \\
\endCD
\tag 3$$
Thus, $\pi_1(\Cal C_i, \eta)^{(\ab)}$ is the geometrically abelian quotient of $\pi_1(\Cal C_i, \eta)$.

\proclaim {Proposition C} The projection  $\pi_1(\Cal C_i, \eta)^{(\ab)} \twoheadrightarrow G_k$ {\bf splits}, $\forall i\in \Bbb N$.
\endproclaim

The author ignores, for the time being, if the projection $\pi_1(\Cal C_{\infty}, \eta)^{(\ab)} \twoheadrightarrow G_k$ splits or not.

\subhead
\S 1. Proof of Theorem A
\endsubhead
In this section we shall prove Theorem A. We use the notations in $\S0$.
We argue by contradiction, and {\bf assume} that the projection
$\pi_1(\Cal C_n,\eta)\twoheadrightarrow G_k$ {\bf splits}, $\forall n\ge 1$.

\proclaim {Proposition 1.1} There exists a relative curve $X\to \Spec \Cal O_k$ with the following properties.

(i) The morphism $X\to \Spec \Cal O_k$ is flat, proper, {\bf stable}, and $X_k\defeq X\times _{\Spec \Cal O_k}\Spec k$ is geometrically connected.

(ii) $X$ is {\bf regular}.

(iii) The set of singular points $X_F^{\text {sing}}$ of the  special fibre $X_F\defeq X\times _{\Spec \Cal O_k}\Spec F$ of $X$ consists
of $F$-{\bf rational ordinary double points}, $U\defeq X_F\setminus X_F^{\text {sing}}$ is $F$-smooth, and 
$U(F)=\emptyset$ holds. 

(iv) $X(k)=\emptyset$ holds.
\endproclaim

\demo{Proof} First, assume $p\neq 2$. Let $\widetilde C\defeq \Bbb P^1_F$ with function field $k(\widetilde C)$. Thus, $\Card (\widetilde C(F))=\Card F+1$ is even.
Arrange the set $\widetilde C(F)$ in pairs of $F$-rational points: $\widetilde C(F)=\{(x_i,y_i)\}_{1\le i\le \frac {\Card F+1}{2}}$. One can identify in $\widetilde C$ the points $x_i$ and $y_i$;
$1\le i\le \frac {\Card F+1}{2}$,
to construct a stable proper $F$-curve $C$ which is geometrically connected and geometrically reduced, with normalisation $\widetilde C\to C$. Moreover, 
the set of singular points $C^{\sing}=\{c_i\}_{1\le i\le \frac {\Card F+1}{2}}$ consists of $F$-rational ordinary double points,
and the pre-image of $c_i$ in $\widetilde C$ consists of the two $F$-rational points $\{x_i,y_i\}$. In particular, $C(F)=C^{\sing}=\{c_i\}_{1\le i\le \frac {\Card F+1}{2}}$.
More precisely, for $1\le i\le \frac {\Card F+1}{2}$, let 
$\widetilde {\Cal O}_i\defeq \Cal O_{\widetilde C,x_i}\cap \Cal O_{\widetilde C,y_i}\subset k(\widetilde C)$, $\Cal N_{x_i}\defeq \frak m_{x_i}\cap \widetilde {\Cal O}_i$, and
$\Cal N_{y_i}\defeq \frak m_{y_i}\cap \widetilde {\Cal O}_i$, where $\frak m_{x_i}$ (resp. $\frak m_{y_i}$) is the maximal ideal of $\Cal O_{\widetilde C,x_i}$ 
(resp. $\Cal O_{\widetilde C,y_i}$). Define $\Cal O_{c_i}\defeq F+\Cal N_{x_i}\Cal N_{y_i}\subset \widetilde {\Cal O}_i$. Then $\Cal O_{c_i}$ is a local ring (with maximal ideal
$\Cal N_{x_i}\Cal N_{y_i}$, and residue field $F$) whose integral closure is $\Cal O_{c_i}$ (cf. [Aubry-Lezzi], Proposition 3.1, Theorem 3.4, and the references therein for  the properties of $\Cal O_{c_i}$, as well as the existence of $C$ with the required properties). 

In case $p=2$. Consider the affine $F$-curve $\Spec (\frac {F[s,t]}{(st)})$, and $\tilde C$ its smooth compactification. Thus, $\tilde C$ consists of two $F$-smooth irreducible components 
$\tilde C_1=\Bbb P^1_F$, and $\tilde C_2=\Bbb P^1_F$, which intersect at the $F$-rational ordinary double point $c=(s,t)\in \Spec (\frac {F[s,t]}{(st)})$.
On each irreducible component $\tilde C_i$ of $\tilde C$; $1\le i\le 2$, the set of $F$-rational points of $\tilde C_i\setminus \{c\}$ is non-empty and comes into pairs of rational points
$\{(x_{i,j},y_{i,j})\}_{1\le j\le \frac {\Card F}{2}}$. 
As above we can identify each of those pairs of $F$-rational points $(x_{i,j},y_{i,j})$ into an $F$-rational ordinary double point $c_{i,j}$
to construct a reducible and geometrically connected stable curve $F$-curve $C$ such that the set of singular points $C^{\sing}$ consists of $F$-rational ordinary double points, a double point $c_{i,j}$ lies on a unique irreducible component of $C$, and $C^{\sing}=C(F)$ (the local ring at $c_i$ is defined as above; the case $p\neq 2$. See. [Rosenlicht], $\S4$, for a discussion of this procedure and the existence of such a curve $C$ 
in the case of reducible curves).

Now the stable $F$-curve $C$ can be deformed to a semi-stable $\Cal O_k$-curve $X\to \Spec \Cal O_k$ with special fibre $X_F=C$ satisfying (i) and (ii)
(cf. [Talpo-Vistoli], Proposition 7.10, Corollary 7.11 and its proof). By our construction (iii) holds also. If $x\in X(k)$, then $X$  specialises in a point $\bar x\in C(F)$
which is a regular point of $C$ and lies on a unique irreducible component of $C$ (cf. [Liu], Corollary 9.1.32). Thus, (iv) follows from (iii).
\qed
\enddemo

Let $X\to \Spec \Cal O_K$ be a regular, proper, flat, and stable $\Cal O_k$-curve as in Proposition 1.1.
Let $y\in X_F(F)$ be an $F$-rational point, which is an ordinary double point  and a regular point of $X$ (cf. Proposition 1.1 (ii) and (iii)).
We fix an isomorphism $\rho:\widehat 
{\Cal O}_{X,y}\isom R[[S,T]]/(ST-\pi)$, and identify $\Cal X\defeq \Spec (\widehat {\Cal O}_{X,y}\otimes _{\Cal O_k}k)$ with $\Cal C_1$ via the isomorphism $\rho_k : \widehat 
{\Cal O}_{X,y}\otimes _{\Cal O_k}k\isom \frac {R[[S,T]]}{(ST-\pi)}\otimes _{\Cal O_k}k$ induced by $\rho$.
Thus, we have scheme morphisms
$$\Cal C_{\infty}\to \cdots \to \Cal C_{n+1}\to \Cal C_n\to \cdots \to \Cal C_1\to X_k.$$

For $n\ge 1$, write 
$$\pi_1(X_k\setminus \Cal S_n,\eta)\defeq \varprojlim _{S_n\subset X_k\setminus \Cal C_n}\pi_1(X_k\setminus S_n,\eta),$$ 
where the projective limit is over all finite sets of closed points $S_n \subset X_k\setminus \Cal C_n$, and $\pi_1(X_k\setminus S_n,\eta)$
is the arithmetic fundamental group of the affine curve $X_k\setminus S_n$ with base point $\eta$. 
(Here, we identify the set of closed points of $\Cal C_n$; $n\ge 1$, with a subset of the set of closed points of $X_k$ which specialise in $y$.)
There is a natural projection
$\pi_1(X_k\setminus \Cal S_n,\eta)\twoheadrightarrow G_k$, and we have a commutative diagram
$$
\CD
\pi_1(\Cal C_n,\eta) @>>> G_k\\
@VVV @| \\
\pi_1(X_k\setminus \Cal S_n,\eta) @>>> G_k\\
\endCD
\tag 4
$$
where the left vertical map is induced by the morphism $\Cal C_n\to  X_k$.

Further, we have a natural map
$$\varprojlim _{n\ge 1} \pi_1(\Cal C_n,\eta)\to \varprojlim _{n\ge 1} \pi_1(X_k\setminus \Cal S_n,\eta),$$ 
and $\varprojlim _{n\ge 1} \pi_1(X_k\setminus \Cal S_n,\eta)$ is naturally identified with the absolute Galois group $G_{k(X)}\defeq \Gal \left ({k(X)}^{\sep}/k(X)\right )$, where
${k(X)}^{\sep}$ is the separable closure of the function field $k(X)$ of $X$ determined by the geometric point $\eta$.

\proclaim {Lemma 1.2}  The projection $G_{k(X)} \twoheadrightarrow G_k$ {\bf splits}.
\endproclaim

\demo{Proof} First, our assumption that the projection $\pi_1(\Cal C_n,\eta)\twoheadrightarrow G_k$ {\bf splits}, implies that the projection 
$\pi_1(X_k\setminus \Cal S_n,\eta) \to G_k$ splits, $\forall n\ge 1$ (cf. diagram (4)).

Let $(H_i)_{i\in I}$ be a projective system of quotients $G_{k(X)}
\twoheadrightarrow H_i$, where $H_i$ sits in an exact sequence $1\to F_i\to H_i\to G_k\to 1$ with $F_i$ finite, and $G_{k(X)}=\varprojlim _{i\in I}H_i$. 
[More precisely, write $G_{k(X)}$ as a projective limit of finite groups $\{\tilde H_i\}_{i\in I}$. Then $\tilde H_i$ fits in an exact sequence $1\to F_i\to \tilde H_i\to G_i\to 1$, where $G_i$ is a quotient 
of $G_k$, and $F_i$ a quotient of $\Gal (k(X)^{\sep}/k(X)\overline k)$. Let $1\to F_i\to H_i\to G_k\to 1$ be the pull-back of the group extension $1\to F_i\to \tilde H_i\twoheadrightarrow G_i\to 1$ 
by $G_k\twoheadrightarrow G_i$. Then $G_{k(X)}=\varprojlim _{i\in I}H_i$]. 
The set $\Sect(G_k,G_{k(X)})$ of group-theoretic sections of the projection $G_{k(X)}\twoheadrightarrow G_k$ is naturally identified with the projective limit $\varprojlim _{i\in I}\Sect(G_k,H_i)$ 
of the sets $\Sect(G_k,H_i)$ of group-theoretic sections of the projection $H_i\twoheadrightarrow G_k$. For each $i\in I$, the set $\Sect(G_k,H_i)$ is non-empty. Indeed,
$H_i$ (being a quotient of $G_{k(X)}$) is a quotient of $\pi_1(X_k\setminus \Cal S_n,\eta)$ for some $n\ge 1$, this quotient $\pi_1(X_k\setminus \Cal S_n,\eta)\twoheadrightarrow H_i$ commutes with the projections onto $G_k$, and we know the projection $\pi_1(X_k\setminus \Cal S_n,\eta) \to G_k$ splits. 
Hence the projection $H_i\twoheadrightarrow G_k$ splits.
Moreover, the set $\Sect (G_k,H_i)$ is, up to conjugation by the elements of $F_i$, a torsor under the group $H^1(G_k,F_i)$ which is finite since $k$ is a $p$-adic local field (cf. [Neukirch-Schmidt-Winberg], (7.1.8) Theorem (iii)). Thus, 
 $\Sect (G_k,H_i)$ is a nonempty finite set. Hence the set $\Sect (G_k,G_{k(X)})$ is nonempty being the projective limit of nonempty finite sets. This finishes the proof of Lemma 1.1.
(See also the proof of Proposition 1.5 in [Sa\"\i di] for similar arguments in a slightly different context.)
\qed
\enddemo

Let $s:G_k\to  G_{k(X)}$ be a section of the projection $G_{k(X)} \twoheadrightarrow  G_k$ (cf. Lemma 1.2).

\proclaim {Lemma 1.3}  The section $s$ is {\bf geometric}, i.e., $s(G_k)\subset D_x$, where $D_x\subset G_{k(X)}$ is a decomposition group associated to a
(unique) rational point $x\in X(k)$. In particular, $X(k)\neq \emptyset$.
\endproclaim

\demo{Proof} This follows from [Koenigsmann], Proposition 2.4 (2).
\qed
\enddemo

 Now the conclusion of Lemma 1.3 that $X(k)\neq \emptyset$ contradicts the assertion (iv) in Proposition 1.1 that $X(k)=\emptyset$. This is a contradiction. Thus, our assumption that 
 the projection $\pi_1(\Cal C_n,\eta)\twoheadrightarrow G_k$ {\bf splits}, $\forall n\ge 1$, can not hold. Let $N\ge 1$ be such that the projection
 $\pi_1(\Cal C_{N},\eta)\twoheadrightarrow G_k$ doesn't splits. Then the projection $\pi_1(\Cal C_n,\eta)\twoheadrightarrow G_k$ doesn't splits, $\forall n\ge N$ as required.
 Indeed, this follows from the fact that for $n\ge N$ we have a natural homomorphism $\pi_1(\Cal C_n,\eta)\to \pi_1(\Cal C_N,\eta)$ which commutes with the projections onto $G_k$.
 Hence if the projection $\pi_1(\Cal C_n,\eta)\twoheadrightarrow G_k$ splits then the projection $\pi_1(\Cal C_N,\eta)\twoheadrightarrow G_k$.
 
This finishes the proof of Theorem A.
\qed

\subhead
\S 2. Proof of Theorem B
\endsubhead
Next, we explain how Theorem B can be derived from Theorem A.
We have, $\forall n\ge 1$, a commutative diagram
$$
\CD
\pi_1(\Cal C_{\infty},\eta) @>>> G_k\\
@VVV @|\\
\pi_1(\Cal C_n,\eta) @>>> G_k\\
\endCD
$$
where the horizontal maps are the natural projections, and the left vertical map is induced by the morphism
$\Cal C_{\infty}\to \Cal C_n$.

Now assume that the projection $\pi_1(\Cal C_{\infty},\eta)\twoheadrightarrow G_k$ splits. Then the projection
$\pi_1(\Cal C_n,\eta)\twoheadrightarrow G_k$ splits, $\forall n\ge 1$, by the above diagram. But this contradicts Theorem A.

This finishes the proof of Theorem B.
\qed

\subhead
\S 3. Proof of Proposition C
\endsubhead
Let $n\ge 1$ be an integer, and $\ell_1$, $\ell_2$, distinct prime integers such that $\ell_1\ge 2n$, and $\ell_2\ge 2n$.
Let $\Cal O_1$, and $\Cal O_2$, be totally ramified extensions of $\Cal O_k$ of degree $\ell_1$, and $\ell_2$, with fraction fields $L_1=\Fr(\Cal O_1)$, and $L_2=\Fr(\Cal O_2)$; respectively.
Thus, the extensions $L_1/k$ and $L_2/k$, are disjoint and $\Cal C_n(L_i)\neq \emptyset$, for $i\in \{1,2\}$. A restriction-corestriction argument shows that the class 
$[\pi_1(\Cal C_n,\eta)^{(\ab)}]$ of the group extension $\pi_1(\Cal C_n,\eta)^{(\ab)}$ in $H^2(G_k,\pi_1(\Cal C_{n,\overline k},\bar \eta)^{\ab})$ is trivial. Thus
the group extension $\pi_1(\Cal C_n,\eta)^{(\ab)}$ splits.

This finishes the proof of Proposition C.
\qed

$$\text{References.}$$

\noindent
[Aubry-Lezzi], Aubry, Y., Lezzi, A., On the Maximum Number of Rational Points on Singular Curves over Finite Fields.
Mosc. Math. J. (2015), Volume 15, Number 4, Pages 615-627, arXiv:1501.03676.

\noindent
[Bourbaki] N. Bourbaki, Alg\`ebre Commutative, Chapitre 9, Masson, 1983.

\noindent
[Grothendieck] Grothendieck, A., Rev\^etements \'etales et groupe fondamental, Lecture 
Notes in Math. 224, Springer, Heidelberg, 1971.

\noindent
[Koenigsmann] Koenigsmann, J., On the section conjecture in anabelian geometry. J. Reine Angew. Math. 588 (2005), 221-235.

\noindent
[Liu] Liu, Q., Algebraic geometry and arithmetic curves, Oxford graduate texts in mathematics 6. Oxford University Press, 2002.



\noindent
[Rosenlicht] Rosenlicht, M., Equivalence relations on algebraic curves, Annals of Mathematics, Vol. 56, No. 1 (Jul., 1952), pp. 169-191.



\noindent
[Sa\"\i di] Sa\"\i di, M., On the existence of non-geometric sections of arithmetic fundamental groups, 
Mathematische Zeitschrift 277, no. 1-2 (2014), 361-372.

\noindent
[Talpo-Vistoli], Talpo, M., Vistoli, A., Deformation theory from the point of view of fibered categories, Handbook of Moduli (Volume III), Editors: Gavril Farkas and Ian Morrison,
publisher Int. Press, Somerville, MA series Adv. Lect. Math. (ALM), volume 26, pages 281-397, arXiv:1006.0497.


\bigskip

\noindent
Mohamed Sa\"\i di

\noindent
College of Engineering, Mathematics, and Physical Sciences

\noindent
University of Exeter

\noindent
Harrison Building

\noindent
North Park Road

\noindent
EXETER EX4 4QF

\noindent
United Kingdom

\noindent
M.Saidi\@exeter.ac.uk

\end
\enddocument